\documentstyle[12pt]{article}
\begin{document}\thispagestyle{empty}\begin{flushright}
OUT--4102--69\\
q-alg/9709031\\
20 September 1997       \end{flushright}\vspace*{2mm}\begin{center}{\Large\bf
Conjectured enumeration of Vassiliev invariants}\vglue 10mm{\large{\bf
D.~J.~Broadhurst                                $^{1)}$}\vglue 4mm
Physics Department, Open University             \\[3pt]
Milton Keynes MK7 6AA, UK     }\end{center}\vfill\noindent{\bf Abstract}\quad
A rational Ansatz is proposed for the generating function $\sum_{j,k}
\beta_{2j+k,2j}x^j y^k$, where $\beta_{m,u}$ is the number of primitive
chinese character diagrams with $u$ univalent and $2m-u$ trivalent vertices.
For $P_m:=\sum_{u\ge2}\beta_{m,u}$, the conjecture leads to the sequence
$$1,1,1,2,3,5,8,12,18,27,39,55,\underline{78,108,150,207,284,388,532,726}$$
for primitive chord diagrams of degrees $m\le20$, with predictions
underlined. The asymptotic behaviour $\lim_{m\to\infty}P_m/r^m=
1.06260548918755$ results, with $r=1.38027756909761$ solving $r^4=r^3+1$.
Vassiliev invariants of knots are then enumerated by $$0,1,1,3,4,9,14,27,44,
80,132,232,\underline{384,659,1095,1851,3065,5128,8461,14031}$$ and
Vassiliev invariants of framed knots by $$1,2,3,6,10,19,33,60,104,184,316,
548,\underline{932,1591,2686,4537,7602,12730,21191,35222}$$ These conjectures
are motivated by successful enumerations of irreducible Euler sums.
Predictions for $\beta_{15,10}$, $\beta_{16,12}$ and $\beta_{19,16}$ suggest
that the action of sl and osp Lie algebras, on baguette diagrams with ladder
insertions, fails to detect an invariant in each case.
\vfill\footnoterule\noindent
$^1$) D.Broadhurst@open.ac.uk;
http://physics.open.ac.uk/$\;\widetilde{}$dbroadhu
\newpage\setcounter{page}{1}

\section{Introduction}

The purpose of this note is to report a conjectured
enumeration of Vassiliev invariants~\cite{Vas}, which results
from fits to the data of~\cite{DBN,JK,JKt} on the enumeration of
primitive chinese character diagrams.
It is motivated by recent successes~\cite{DJB,CRM}
of conjectured enumerations~\cite{BK14} of irreducible multiple zeta values
(MZVs)~\cite{DZ,BBG,BG} and their extensions, alternating Euler
sums~\cite{Eul,BBB}. It is consistent with all 64 proven results~\cite{JK}
and lower bounds~\cite{JKt} on the bigrading
of chinese characters by degree, $m\le14$,
and number of univalent vertices, $u=2j\le m$.
Moreover, it saturates the bounds for $m=13,14$.
For $m\ge15$, some lower bounds are exceeded by the
Ansatz, suggesting invariants that were not detected by
Lie algebras of type sl and osp operating on a restricted set of
diagrams.

A chinese character diagram of degree $m$, with $u$ univalent vertices,
has $2m-u$ trivalent vertices. If the univalent vertices
are attached to a circle, representing a closed Wilson line,
the resultant graph has $m+1$ loops,
when regarded as a Feynman diagram in the perturbative
expansion of a topological quantum field theory~\cite{EW,DBNth}.
In abstracting from the notion of a topological field theory,
Dror Bar--Natan was led to the consideration of weight systems~\cite{DBNth}.
Combinatorical investigation of such weight systems enables
an enumeration of Vassiliev invariants for knots, via the equivalence
proven in~\cite{DBN,MK,DBNS}.

The data of~\cite{DBN,JK} concern the number, $\beta_{m,u}$,
of diagrams that cannot be reduced merely by the application of the
antisymmetry and Jacobi relations (called AS and IHX relations in~\cite{DBN})
that obtain when structure constants of a Lie algebra are associated
with the vertices of a chinese character.
The numbers so obtained classify the dimensionality
of the naturally graded spaces of primitive elements in the bialgebra
established by weight systems~\cite{DBN}. Summing over $u\ge2$,
one obtains an enumeration of connected chord diagrams, modulo four-term
relations, which is equivalent to
an enumeration of primitive Vassiliev invariants by~\cite{DBN}
\begin{equation}
P_m:=\sum_{u\ge2}\beta_{m,u}\,.\label{pm}
\end{equation}
These, in turn, generate the numbers, $V_m$ and $F_m$, of Vassiliev
invariants of knots and framed knots, respectively, via the Euler transforms
\begin{eqnarray}
\prod_{m\ge2}(1-y^m)^{-P_m}&=&1+\sum_{m\ge2}V_my^m\,,\label{vm}\\
\prod_{m\ge1}(1-y^m)^{-P_m}&=&1+\sum_{m\ge1}F_my^m\,,\label{fm}
\end{eqnarray}
with $P_1=\beta_{1,2}=1$ implying that $V_{m+1}=F_{m+1}-F_m$.

It is known~\cite{Vog} that at some degree, $m>12$,
one will eventually encounter structure more general than that spanned
by Lie algebras. Recently~\cite{JL} it was shown that a Lie superalgebra
detects structure at $m=19$ and $u=4$
that is not detected by semisimple Lie algebras.
Moreover, even superalgebras are proven
to be insufficient, in the long run~\cite{Vog}.

As a contribution to such studies, we here attempt to infer,
from existing data~\cite{JK,JKt}, a closed rational form for the generator
\begin{equation}
b(x,y):=\sum_{j,k\ge0}(\beta_{2j+k,2j}-1)x^j y^k\,.\label{bxy}
\end{equation}
Assuming that Vassiliev invariants
fail to distinguish a knot from its mirror image we set
$\beta_{m,2j+1}=0$. Then, with $P_1=\beta_{1,2}=\beta_{0,0}=1$,
we conjecture the sequences generated by~(\ref{pm},\ref{vm},\ref{fm}).
The key ingredient is the hypothesis of common features
in the generators for primitive Vassiliev invariants and
for irreducible MZVs and Euler sums.

Like the chinese character diagrams of~\cite{DBN},
the MZVs of~\cite{DZ} and the Euler sums of~\cite{Eul}
have a bigrading.
In the case of MZVs and Euler sums, it is by weight $w$ and depth $d$.
The distinctive prediction
of the conjecture of~\cite{BK14} is that the number, $D_{w,d}$,
of irreducible MZVs is generated by a
pseudopolynomial at any fixed depth. In other words,
the generating function $\sum_{j}D_{2j+3d,d}x^j$
is believed to have singularities
only at roots of unity, for fixed $d$.

The clue that the enumerations of MZVs and Vassiliev invariants might be
akin comes from the Kontsevich integral~\cite{MK},
which is a universal knot invariant, corresponding
to the holonomy of the flat Knizhnik-Zamolodchikov connection~\cite{KZ,Drin}.
It thus evaluates in terms of MZVs, in a natural manner \cite{CK}.
Moreover Dirk Kreimer, the present author and coworkers
have accumulated abundant evidence~\cite{DK1}--\cite{4TR} that
knots and MZVs are connected via the counterterms of quantum field theory.
These counterterms have recently
been construed~\cite{DKh} as deriving from the antipode of
a quasi-Hopf algebra, whose failure to be coassociative
may be controlled by a Drinfeld associator~\cite{Drin,CK},
specific to the quantum
field theory in question.
One is thus encouraged to fit the recent~\cite{JK}
bigraded data on chinese characters with pseudopolynomial
generating functions, deriving from an underlying rational generator
in~(\ref{bxy}).
This was surprisingly easy to achieve, and quickly suggested
a closed form for the rational generator.

Section~2 gives a detailed analysis of the data~\cite{DBN,JK}
on $\beta_{m,u}$, leading to pseudopolynomial generators
for fixed values of $m-u$ or $u=2j$.
In Section~3 these observations are subsumed into a relatively simple
Ansatz that fits all available data and satisfies all known bounds.
Section~4 gives comments and conclusions.

\section{Analysis of data}

The impressive progress reported by Jan Kneissler in~\cite{JK}
extended the results on $\beta_{m,u}$ from degrees $m\le9$,
achieved in~\cite{DBN},
up to $m\le12$. All values of $\beta_{m,u}$ with
$m\le12$ are known, save $\beta_{12,0}$, which is merely bounded
from below.
In addition to these, one knows that~\cite{DBN}
\begin{eqnarray}
\beta_{m,0}&=&\beta_{m+1,2}\,,\label{u2}\\
\beta_{2j,2j}&=&1\,,\label{k0}\\
\beta_{2j+1,2j}&=&\left\lfloor\frac{j+3}{3}\right\rfloor\,,\label{k1}
\end{eqnarray}
where $\lfloor x \rfloor$ denotes the largest integer that does not
exceed $x$.

Moreover, there exist lower bounds, obtained
by a `thickening' procedure~\cite{DBN} that enumerates the
invariants detected by sl and osp Lie algebras. Restricting
attention to diagrams generated by ladder insertions in baguette~\cite{CD}
diagrams (called caterpillar diagrams in~\cite{JK}) lower
bounds on $\beta_{m,u}$ were obtained in~\cite{JKt} for
$20\ge m\ge13$. With $u=2$, further invariants were found
for $20\ge m\ge14$.
The values in Table~1 that are lower bounds
are indicated by underlining. The Ansatz adduced in
Section~3 saturates these.

\begin{table}[htb]\caption{
Values and \underline{lower bounds} for $\beta_{m,u}$ with $m\le14$.}
\[\begin{array}{l|rrrrrrrrr}
&u=0&u=2&u=4&u=6&u=8&u=10&u=12&u=14\\\hline
m=\phantom{1}   0&1                \\
m=\phantom{1}   1&1&1              \\
m=\phantom{1}   2&1&1              \\
m=\phantom{1}   3&1&1              \\
m=\phantom{1}   4&2&1& 1           \\
m=\phantom{1}   5&2&2& 1           \\
m=\phantom{1}   6&3&2& 2& 1        \\
m=\phantom{1}   7&4&3& 3& 2        \\
m=\phantom{1}   8&5&4& 4& 3& 1     \\
m=\phantom{1}   9&6&5& 6& 5& 2     \\
m=             10&8&6& 8& 8& 4&1   \\
m=             11&9&8&10&11& 8&2   \\
m=12&\underline{11}&9&13&15&12&5&1 \\
m=13&\underline{13}&\underline{11}&\underline{16}&\underline{20}&
\underline{18}&\underline{10}&3&\\
m=14&\underline{15}&\underline{13}&\underline{19}&\underline{25}&
\underline{26}&\underline{17}&7&1
\end{array}\]
\end{table}

The possibility of a simple enumeration of Vassiliev invariants
was suggested by observing that Oliver Dasbach's result~\cite{Das}
for the third diagonal of Table~1, namely
\begin{equation}
\beta_{2j+2,2j}=
\left\lfloor\frac{(j+3)^2+3}{12}\right\rfloor\,,\label{k2}
\end{equation}
is echoed by a simple fit to the first two columns of Table~1, namely
\begin{equation}
\beta_{m,0}-1=\beta_{m+1,2}-1\,\stackrel{?}{=}\,\left\lfloor
\frac{(m-1)^2+3}{12}\right\rfloor\,,\label{u0}
\end{equation}
which is proven for $11\ge m\ge0$. This conjecture\footnote{Conjectured
equalities are indicated by $\,\stackrel{?}{=}\,$, as in~(\ref{u0}).} was
lodged as A014591 in Neil Sloane's encyclopedia of integer
sequences~\cite{NJAS}, prior to the lower bounds of~\cite{JKt} for
$\beta_{m,0}$ with $19\ge m\ge12$.
It was then found to saturate these lower bounds for $16\ge m\ge12$.
Having thus fitted 12 data and saturated a further 5 bounds, we suggest
that the lower bounds of~\cite{JKt} for $\beta_{m,0}$
are exceeded for $m=17,18,19$, in accordance with~(\ref{u0}).

Encouraged by results for the three leading diagonals of Table~1,
and fits to its first two columns, we conjectured
that the generator~(\ref{bxy}) has Taylor coefficients in $x$ that
are pseudopolynomials in $y$, with singularities at $y^2=1$
and $y^3=1$, and Taylor coefficients in $y$ that are pseudopolynomials
in $x$, with singularities at $x^2=1$ and $x^3=1$.

The guiding intuition
behind this conjecture was the fact that the conjectured
enumeration of MZVs in~\cite{BK14} had been seeded by singularities
at square and cube roots of unity, with
\begin{equation}
\prod_{j\ge0}\prod_{d>0}(1-x^j y^d)^{D_{2j+3d,d}}\;\stackrel{?}{=}\;
1-\frac{y}{1-x}-\frac{y^2}{1-x^2}\frac{y^2-x^3}{1-x^3}\,,\label{MZV}
\end{equation}
proposed~\cite{BK14} as the generator of the number, $D_{w,d}$, of irreducible
MZVs of weight $w\ge3d$ and depth $d$. This Ansatz was initially
inferred from rather limited data. Subsequently, it has survived intensive
testing. For example the prediction~\cite{BK14} $D_{23,7}=4$
was confirmed~\cite{CRM}
as the rank deficiency of the $447\,678\times74\,613$ matrix that results
from the interplay of the weight-length~\cite{DZ,CK}
and depth-length~\cite{BBB}
shuffle algebras for MZVs. It is an outstanding puzzle to understand
the third term in~(\ref{MZV}), which is absent in the enumeration
\begin{equation}
\prod_{j\ge0}\prod_{d>0}(1-x^j y^d)^{M_{2j+3d,d}}\;\stackrel{?}{=}\;
1-\frac{y}{1-x}\,,\label{Eul}
\end{equation}
for the number, $M_{w,d}$, of irreducible Euler sums of weight $w$
and depth $d$ that serve to reduce all MZVs. The difference
between~(\ref{MZV}) and~(\ref{Eul}) first surfaces at weight $w=12$,
where $\zeta(4,4,2,2):=\sum_{k>l>m>n>0}k^{-4}l^{-4}m^{-2}n^{-2}$ is
{\em not\/} reducible to MZVs of depth $d<4$ but {\em is\/} reducible to
the depth-2 alternating Euler sum~\cite{Eul}
$\sum_{m>n>0}(-1)^{m+n}m^{-9}n^{-3}$,
as anticipated~\cite{BGK} by the connection~\cite{DK1} between knots and
numbers effected by the counterterms of quantum field theory.

Close study of Table~1 suggests that the enumeration of
Vassiliev invariants, with $m-u\ge3$, is {\em not\/} obtainable
from that of irreducible MZVs by mere relabelling.
This presumably reflects a difference in
the Drinfeld associators controlling the failures of coassociativity
in the quasi-Hopf algebras that relate to knot theory,
in the case of rational-valued Vassiliev invariants~\cite{DBNa}, and
to the Knizhnik-Zamolodchikov equation~\cite{KZ}, in the case of
irrational MZVs. However, the connections between knots and MZVs,
effected by the Kontsevich integral~\cite{MK} and by quantum
field theory~\cite{DK1}, lead us to expect resemblances between
generating functions for these distinct structures.

Hence we sought to generate columns and diagonals of Table~1 with
pseudopolynomials whose singularities occur only at square and cube
roots of unity.
With such limited pseudopolynomial singularities posited for
\begin{eqnarray}
g_k(x)&:=&\sum_{j\ge0}\beta_{2j+k,2j} x^j\,,\label{gk}\\
h_j(y)&:=&\sum_{k\ge0}(\beta_{2j+k,2j}-1) y^k\,,\label{hj}
\end{eqnarray}
it is apparent that the full generator~(\ref{bxy}) must have a
denominator that involves a factor which couples $x$ to $y$,
in order to generate terms that are more singular at $x=1$
when one expands to higher order in $y$, and vice versa.

Fortunately, the proven sequence~\cite{JK}
\begin{equation}
1,1,1,2,3,5,8,12,18,27,39,55\label{rat}
\end{equation}
for $P_m$, with $m\le12$, provided a clue as to this
coupling, since the ratio of
successive terms appears to tend to a value $r\approx1.4$.
An origin of such a ratio is not hard to imagine. Consider, for example,
the number, $D_{3d,d}$, of irreducible MZVs of depth $d$ and weight
$w=3d$, which is the lowest weight at which sums of this depth
may exhibit irreducibility. Then from~(\ref{MZV}), at $x=0$, one obtains
\begin{equation}
\prod_{d>0}\frac{1}{(1-y^d)^{D_{3d,d}}}\;\stackrel{?}{=}\;
\frac{1}{1-y-y^4}\,,\label{d3d}
\end{equation}
which has been verified through $O(y^7)$, i.e.\ for weights $w=3d\le21$,
and is lodged as A020999 in~\cite{NJAS}.
The asymptotic growth is determined by the singularity of $1/(1-y-y^4)$
at $y=1/r$, with $r=\lim_{d\to\infty} D_{3d+3,d+1}/D_{3d,d}=1.38027756909761$
solving $1=r^{-1}+r^{-4}$. Moreover, one can generate such
an asymptotic growth in the case of~(\ref{pm}), as follows.

Suppose that $b(x,y)$
in~(\ref{bxy})
had a singularity of the form $1/(1-y-x^2)$.
Then the summation in~(\ref{pm}) would detect a $1/(1-y-y^4)$
singularity from the term with $x=y^2$ in
\begin{equation}
p(y):=\sum_{m\ge1}(P_m-1)y^m=b(y^2,y)-b(0,y)
+\frac{y^4}{(1-y)(1-y^2)}\,,\label{py}
\end{equation}
and hence give $\lim_{m\to\infty}P_{m+1}/P_{m}=r\approx1.38$.
Thus previous experience with MZVs, combined
with a ratio test $r\approx1.4$ from~(\ref{rat}),
led us to suspect a $1/(1-y-x^2)$ singularity in~(\ref{bxy}).
Hence we expected that expanding in powers of $y$  would lead
to higher powers of $1/(1-x^2)$ in~(\ref{gk}), as $k$ increases.
This observation greatly facilitated the discovery of simple
pseudopolynomial fits to the first few instances of~(\ref{gk}), namely
\begin{eqnarray}
g_0(x)&=&\frac{1}{1-x}=\sum_j D_{3+2j,1} x^j\,,\label{g0}\\
g_1(x)&=&\frac{g_0(x)}{1-x^3}=\sum_j D_{8+2j,2} x^j\,,\label{g1}\\
g_2(x)&=&\frac{g_1(x)}{1-x^2}=\sum_j(D_{11+2j,3}-D_{8+2j,2})x^{j-1}\,,
\label{g2}\\
g_3(x)&\;\stackrel{?}{=}\;&\frac{g_2(x)}{1-x^2}
+\frac{x}{(1-x^2)(1-x^3)}\,,\label{g3}\\
g_4(x)&\;\stackrel{?}{=}\;&\frac{g_3(x)}{1-x^2}
+\frac{1}{1-x^3}\,,\label{g4}\\
g_5(x)&\;\stackrel{?}{=}\;&\frac{g_4(x)}{1-x^2}
+\frac{x}{1-x^2}\,,\label{g5}
\end{eqnarray}
for the first 6 diagonals of Table~1.
For the first 4 columns we obtained the fits
\begin{eqnarray}
h_0(y)=y h_1(y)&\;\stackrel{?}{=}\;&\frac{y^4}{(1-y)(1-y^2)(1-y^3)}
\,,\label{h1}\\
y h_2(y)&\;\stackrel{?}{=}\;&(1+y)h_1(y)\,,\label{h2}\\
y h_3(y)&\;\stackrel{?}{=}\;&(1+y^2)h_2(y)\,,\label{h3}
\end{eqnarray}
by the mere device of inserting $1/(1-y^3)$ in the transparent fits
\begin{eqnarray}
\overline{h}_0(y)=y \overline{h}_1(y)&\;\stackrel{?}{=}\;&
\frac{y^4}{(1-y)(1-y^2)}\,,\phantom{(1-y^3)}\label{hb1}\\
y\overline{h}_2(y)&\;\stackrel{?}{=}\;&
(1+y)\overline{h}_1(y)\,,\label{hb2}\\
y\overline{h}_3(y)&\;\stackrel{?}{=}\;&
(1+y^2)\overline{h}_2(y)\,.\label{hb3}
\end{eqnarray}
to {\em all\/} of the results with $m\le20$ and $u\le6$
given in~\cite{JKt} for the contributions to~(\ref{hj}) coming from
orientable surfaces, i.e.\ from Lie algebras of type sl.

The challenge was then to use the posited factor
$1/(1-y-x^2)$ to marry~(\ref{g0}--\ref{g5}) with~(\ref{h1}--\ref{h3}).
This was achieved as follows.

\section{Conjecture}

It is conjectured that the numbers of primitive Vassiliev invariants
are generated by
\begin{equation}
b(x,y):=\sum_{j,k\ge0}(\beta_{2j+k,2j}-1)x^j y^k
\;\stackrel{?}{=}\;\frac{b_0y^4+b_1xy^3+b_2x^2y^2}{1-x^3}
+\frac{b_3x^3y+b_4x^4}{(1-x^3)(1-y-x^2)}\,,\label{B}
\end{equation}
with
\begin{equation}
b_0=b_1=\frac{b_2}{1+y}=\frac{b_3}{1-y^3}=1+b_4=
\frac{1}{(1-y)(1-y^2)(1-y^3)}\,.\label{b1}
\end{equation}
This is consistent with all currently available values and bounds for
$\beta_{m,u}$ and leads to
\begin{equation}
p(y):=\sum_{m\ge1}(P_m-1)y^m\;\stackrel{?}{=}\;
\frac{y^4-y^8-y^{10}-y^{12}-y^{17}}
{(1-y)(1-y^2)(1-y^3)(1-y^6)(1-y-y^4)}\,,\label{pis}
\end{equation}
with $P_{m}:=\sum_{u\ge2}\beta_{m,u}$
enumerating primitive chord diagrams.

\section{Comments and conclusions}

We note the following features of~(\ref{B},\ref{b1}) and their
corollary~(\ref{pis}).
\begin{enumerate}
\item It was assumed that $(1-y)(1-y^2)(1-y^3)(1-x^3)(1-y-x^2)b(x,y)$
is polynomial in $x$ and $y$ and that
$b(\infty,y)=b(x,\infty)=0$. This was motivated by the
MZV analogies in~(\ref{g0}--\ref{g2}) and by
the pseudopolynomial fits of~(\ref{g3}--\ref{h3}).
\item Thus Ansatz~(\ref{B},\ref{b1}) was fully determined by the data
of Table~1 with $m-6\le u\le8$. The remaining 29 data of Table~1
were successful predictions. It is remarkable that the
simple pseudopolynomials of~(\ref{b1}) fit the entirety of
Table~1 and satisfy all further results and bounds.
If the generator proves to be
somewhat different, it will be interesting to see how it
preserves the existing success, with comparable economy.
\item Summing over $u\ge2$ we obtain the tallies
\begin{eqnarray*}
P_{13}\;\stackrel{?}{=}\;11+16+20+\phantom{1}18+\phantom{1}10+\phantom{10}3
\phantom{+106+52+12+1}\,&=&\phantom{1}                                   78\\
P_{14}\;\stackrel{?}{=}\;13+19+25+\phantom{1}26+\phantom{1}17+\phantom{10}7
+\phantom{10}1\phantom{+52+12+1}\,&=&                                   108\\
P_{15}\;\stackrel{?}{=}\;15+23+31+\phantom{1}35+\phantom{1}28+\phantom{1}15
+\phantom{10}3\phantom{+52+12+1}\,&=&                                   150\\
P_{16}\;\stackrel{?}{=}\;17+27+38+\phantom{1}46+\phantom{1}42+\phantom{1}28
+\phantom{10}8+\phantom{1}1\phantom{+12+1}\,&=&                         207\\
P_{17}\;\stackrel{?}{=}\;20+31+45+\phantom{1}60+\phantom{1}60+\phantom{1}46
+\phantom{1}19+\phantom{1}3\phantom{+12+1}\,&=&                         284\\
P_{18}\;\stackrel{?}{=}\;22+36+53+\phantom{1}75+\phantom{1}83+\phantom{1}72
+\phantom{1}36+10+\phantom{1}1\phantom{+1}\,&=&                         388\\
P_{19}\;\stackrel{?}{=}\;25+41+62+\phantom{1}93+          111+          107
+\phantom{1}64+25+\phantom{1}4\phantom{+1}\,&=&                         532\\
P_{20}\;\stackrel{?}{=}\;28+46+71+          114+          144+          152
+106+52+12+1\!&=&                                                       726
\end{eqnarray*}
which generate the unproven parts of the
sequences in the abstract.
\item The parallel between orientable lower bounds,
generated by~(\ref{hb1},\ref{hb2},\ref{hb3}) for $m\le20$,
and the fits of~(\ref{h1},\ref{h2},\ref{h3}) to Table~1, with $m\le14$,
is quite remarkable.
\item The lower bounds of~\cite{JKt} are saturated at $m\le14$,
for all $u$, and at $m\le17$, for $u=2$ and for $u=m-3$.
\item With $m=15$, the fit yields values of $\beta_{m,u}$
that exceed, by unity, the lower bounds of~\cite{JKt} for $u=4,6,8,10$,
from invariants detected by Lie algebras of type sl and osp, restricted
to baguette diagrams with ladder insertions.
For $u=2$, it is known that such bounds are exceeded at all $m\ge14$.
\item The expectations that
$\beta_{19,16}=25$, $\beta_{16,12}=28$, $\beta_{15,10}=28$
are the first cases where~(\ref{g3},\ref{g4},\ref{g5})
require an invariant additional to the lower bounds of~\cite{JKt}.
It would be very interesting to know what {\em upper\/} bound
algorithms~\cite{JK} yield in these cases.
\item The numerator of~(\ref{pis}) is easy to characterize:
the coefficients of $y^m$ are determined by proven
results for $m\le12$; thereafter they vanish, save at
$m=17$, where a contribution results from the final
term in~(\ref{py}), which dominates at large $y$.
This gives $\lim_{y\to\infty}p(y)/y=1$, coming
from the sole contribution with $m<u$, namely $\beta_{1,2}=1$,
which mimics the appearance of $\pi^2$ in the reduction of MZVs.
\item Were the predictions for $P_m$ to
be confirmed for $17\ge m\ge13$,
those for $m\ge18$ would appear to be rather secure,
lending weight to the asymptotic prediction
\begin{equation}
\lim_{m\to\infty}P_m/r^m=
\lim_{y\to1/r}(1-r y)p(y)\;\stackrel{?}{=}\;1.06260548918755\ldots\label{C}
\end{equation}
with $r=1.38027756909761\ldots$ solving $r^4=r^3+1$. This comes from
the singularity at $y=1/r$ of the $1/(1-y-y^4)$ term in~(\ref{pis}),
which results from setting $x=y^2$ in~(\ref{B}), in accordance
with~(\ref{py}). The analogous singularity in~(\ref{d3d})
gives $\lim_{d\to\infty}D_{3d+3,d+1}/D_{3d,d}=r$.
\item The rational fit to Table~1 suggests that so-called {\em primitive\/}
Vassiliev invariants, enumerated by $\beta_{m,u}$,
are analogous to the elements in a search space
for reducing MZVs to {\bf Q}-linear combinations
of basis terms. In the case of MZVs, the basis terms enumerated
by a rational generator include both
irreducible MZVs, such as $\zeta(5,3):=\sum_{m>n>0} m^{-5} n^{-3}$,
and their {\em products}, such as $\zeta(5)\zeta(3)$.
If the successful phenomenology
of the present enterprise be more than
accidental, then there might be a further decomposition
of the invariants enumerated by $\beta_{m,u}$,
so far undetected in the rational domain of Vassiliev
invariants, while its counterpart
is readily detectable by the PSLQ algorithm~\cite{CRM,PSLQ}
in the irrational domain of MZVs.
\end{enumerate}
{\bf Acknowledgements}: I thank David Bailey, Dror Bar--Natan,
Jonathan Borwein, Pierre Deligne, Pavel Etingof, Sasha Goncharov,
John Gracey, Dirk Kreimer, Petr Lisonek, Justin Roberts, Neil Sloane,
Ioannis Tsohantjis, Simon Willerton and Don Zagier, for generous advice.
Most importantly, Dirk Kreimer encouraged me to believe that enumerations of
Vassiliev invariants and irreducible MZVs might be akin.
Jan Kneissler's website data provided the vital stimulus.
This work was completed at the 18th UK High Energy Physics Institute,
in St Andrews, thanks to discussions with Dirk Kreimer,
enabled by HUCAM grant CHRX--CT94--0579.

\newpage
\raggedright


\begin{thebibliography}{99}
\bibitem{JL} J.\ Lieberum,
{\sl On Vassiliev invariants not coming from semisimple Lie algebras},
{\tt q-alg/9706005}.
\end{thebibliography}
\end{document}